\documentclass[11pt]{article}
\usepackage{amssymb}
\newtheorem{lem}{Lemma}[section]

\newtheorem{pro}{Proposition}[section]
\newtheorem{thm}{Theorem}[section]

\makeatletter\@addtoreset{equation}{section}
\renewcommand\theequation{\thesection.\@arabic\c@equation}

\newenvironment{Literature}[1]
{\begin{thebibliography}{#1}}{\end{thebibliography}}

\begin{document}

\begin{center}
{\LARGE \bf  Remark on the Limit Case of  Positive Mass }

\bigskip
{\LARGE \bf Theorem  for Manifolds with Inner Boundary}

\bigskip  \bigskip
{\large  Eui Chul Kim}
\end{center}

\bigskip   \noindent
Department of Mathematics, Inha University  \\ Yonghyoundong 253,
402-751 Inchon,  Korea    \\ e-mail: eckim@chol.com

\bigskip \bigskip \noindent
{\bf Abstract:}  In [5] Herzlich proved a new positive mass
theorem for Riemannian 3-manifolds $(N, g)$ whose mean curvature
of the boundary allows some positivity. In this paper we study what happens to
 the limit case of the theorem when, at a point of the boundary,  the smallest
positive eigenvalue of the Dirac operator of the boundary is
strictly larger than one-half of the mean curvature (in this case the mass $m(g)$ must be
strictly positive). We prove that the mass is bounded from below
by a positive constant $c(g),  \, m(g)  \geq  c(g)$, and
the equality $m(g) = c(g)$ holds only if, outside a compact  set, $(N,
g)$ is conformally flat and the scalar curvature vanishes. The
constant $c(g)$ is uniquely determined by the metric $g$ via a
Dirac-harmonic spinor.

\bigskip \noindent                                                                                                                                                                 {\bf MSC(2000):}  53C27, 83C40    \\
{\bf Keywords:}  Positive mass theorem, Harmonic spinor, Rigidity
\\

\bigskip  \bigskip
\noindent
\section{Introduction}

Let $(N, g)$ be a complete Riemannian
3-manifold with boundary which is diffeomorphic to the Euclidean space ${\mathbb R}^3$ minus an
open 3-ball centered at the origin.   Let $ r(y) = \sqrt{\sum_{i=1}^3 y_i^2},   \, y = (y_1,  y_2, y_3  ) \in
{\mathbb R}^3$, be the standard distance function to the origin of ${\mathbb
R}^3$.  Then $(N, g)$ is
called {\it asymptotically flat} of order $\tau >   \frac{1} {2}$, if
 there is a diffeomorphism $ \Phi : N  \longrightarrow {\mathbb R}^3 \backslash \mbox{\{an open 3-ball\}}$ \, such that the coefficients of the metric $g$ in the induced rectangular coordinates satisfy
\[       g_{ij} = \delta_{ij} + O(r^{-\tau}),  \qquad   g_{ij,k} =
O(r^{-\tau-1})  , \qquad   g_{ij,k,l} = O( r^{-\tau-2})  
\]
as $r =  r(\Phi)  \longrightarrow \infty$.
Let $S(r)  \subset  N$ denote the $\Phi$-inverse image of a round 2-sphere in ${\mathbb R}^3$, centered at the origin and of sufficiently large radius $r>0$.
Throughout the paper we identify
\[    N  =  \bigcup_{r \geq r_o}  S(r)   \quad    \mbox{for some fixed constant} \ r_o > 0   .   \]
The {\it mass} of
$(N, g)$ is usually defined  by [1]        \begin{equation}        m(g)   =   \frac{1}{16 \pi} \, \lim_{r
\to \infty}   \sum_{i,j=1}^3  \int _{S(r)} ( g_{ij,j} -  g_{jj,i} ) \nu^i dS,         \end{equation}
where $\nu$ is the outward unit normal to spheres $S(r)$ and  $dS$ is the
 area form of spheres $S(r) \subset N$.   We remark here that one can express this
definition in a coordinate-independent way, by considering a flat metric on $N$ as a reference metric. Namely, let $g_{\rm eu}$ be a metric on $N$ which is the
pullback of the Euclidean metric on ${\mathbb R}^3 \backslash \mbox{\{an open 3-ball\}}$ via the diffeomorphism
$ \Phi : N  \longrightarrow {\mathbb R}^3 \backslash \mbox{\{an open 3-ball\}}$.  Then the equation (1.1) is in fact equal to
 \begin{eqnarray}
 m(g)  &  =  &  \frac{1}{16 \pi}   \,   \lim_{r \to \infty} \int _{S(r)}
g_{\rm eu} ( {\rm div}_{g_{\rm eu}} (g) - {\rm grad}_{g_{\rm eu}} ( {\rm Tr}_{g_{\rm eu}} (g) ), \, V_{\rm eu} ) \mu_{S(r)}(g_{\rm eu}),      \\
&          &        \nonumber      \\
& = & \frac{1}{16 \pi}   \,   \lim_{r \to \infty} \int _{S(r)}
g ( {\rm div}_{g_{\rm eu}} (g) - {\rm grad}_{g_{\rm eu}} ( {\rm Tr}_{g_{\rm eu}} (g) ), \, V_g ) \mu_{S(r)}(g),
\end{eqnarray}
where $V_{\rm eu}$ (resp. $V_g$) is the outward unit normal to spheres $(S(r), g_{\rm eu})$ (resp.  $(S(r), g)$) and  $\mu_{S(r)}(g_{\rm eu})$
(resp. $\mu_{S(r)}(g)$) is the area form of spheres $(S(r), g_{\rm eu})$ (resp.  $(S(r), g)$). When one applies the Witten-type spinor method to prove
positivity of the mass, one should use the latter equation (1.3) [2, 5, 6, 9, 11].
Note that the equations (1.2)-(1.3)  are independent of deformation of the foliation  
 $N  =  \bigcup_{r \geq r_o}  S(r) $ via a diffeomorphism $F : N  \longrightarrow  N$  , since Stokes' theorem implies that
\begin{eqnarray*}
m(g)  &   =    &     \frac{1}{16 \pi}   \,  \lim_{r \to \infty} \int _{S(r)}
g ( {\rm div}_{g_{\rm eu}} (g) - {\rm grad}_{g_{\rm eu}} ( {\rm Tr}_{g_{\rm eu}} (g) ), \, V_g ) \mu_{S(r)}(g)     \\
&       &      \\
& = &   \frac{1}{16 \pi}    \,  \int _{\partial N}
g ( {\rm div}_{g_{\rm eu}} (g) - {\rm grad}_{g_{\rm eu}} ( {\rm Tr}_{g_{\rm eu}} (g) ), \, V_g ) \mu_{S(r)}(g)     \\
&        &       \\
&        &   +     \frac{1}{16 \pi}     \,   \int _N   {\rm div}_g    \Big\{
 {\rm div}_{g_{\rm eu}} (g) - {\rm grad}_{g_{\rm eu}} ( {\rm Tr}_{g_{\rm eu}} (g) )   \Big\}    \mu_{S(r)}(g)   
\end{eqnarray*} 
 whose right-hand side is independent of a choice of foliation on $N$ by 2-spheres.

\bigskip
The mass is a geometric invariant of Riemannian asymptotically
flat manifolds and of importance in Riemannian geometry  as well
as in general relativity.  In [3, 7] one finds  an excellent exposition of the positive mass conjecture as well as the Penrose conjecture and a full list
of related papers.
A fundamental problem about the mass is
to investigate the relation between the scalar curvature $S_g$ of
the manifold $(N, g)$, the mean curvature ${\rm Tr}_g(\Theta)$ of the inner
boundary $(\partial N, g \big\vert_{\partial N} )$ and the mass
$m(g)$ (Here $\Theta$ indicates the second fundamental form of the boundary). The Riemannian positive mass theorem, proved by Schoen and
Yau [10], states that, if $(N, g)$ is an asymptotically flat
3-manifold of non-negative scalar curvature $S_g \geq 0$ with
minimal boundary ${\rm Tr}_g(\Theta)  \equiv  0$, then the mass is
non-negative $m(g)  \geq 0$. In fact, the limit case of zero mass
can not be attained and so the mass must be strictly positive. The
Penrose conjecture, recently proved by Huisken and Ilmanen [7],
improves the positive mass theorem and states that, if the
boundary is not only minimal but also outermost (i.e., $N$
contains no other compact minimal  hypersurfaces), then
\[   m(g)   \  \geq  \  4  \sqrt{\frac{{\rm Area}(\partial N, g)}{\pi} }       \]
with equality if and only if $(N, g)$ is isometric
to the spatial Schwarzschild manifold.

\bigskip
In [5] Herzlich proved a new positive mass theorem for manifolds
with inner boundary (see Theorem 2.1), making use of
Dirac-harmonic spinors with well-chosen spectral boundary
condition (see the PDE system (2.7) below).    A
remarkable feature of the theorem is that the mass $m(g)$ is
non-negative even if there is some positivity of the mean
curvature of the boundary. The limit case of zero mass (the flat
space) occurs only if the smallest positive eigenvalue $\lambda$
of the Dirac operator of the boundary is equal to  one-half of the
mean curvature ${\rm Tr}_g (\Theta)$, i.e.,
\[      \lambda   =  2 \,  \sqrt{ \frac{\pi}{ {\rm Area} (\partial N, g) } } =  \frac{1}{2}  \, {\rm Tr}_g (\Theta)  .  \]

\par
The object of this paper is to study what happens to the limit case of the theorem when
\[      2 \,  \sqrt{ \frac{\pi}{ {\rm Area} (\partial N, g) } }  \,  \geq   \,    \frac{1}{2}  \,  \sup_{\partial N}  \{ {\rm Tr}_g (\Theta) \}
 \qquad     \mbox{and}    \qquad    2 \,  \sqrt{ \frac{\pi}{ {\rm Area} (\partial N, g) } }   \not\equiv  \frac{1}{2}  \, {\rm Tr}_g (\Theta)  ,   \]
in which case the zero mass $m(g) = 0$ can not be attained. We
will prove (see Theorem 3.1) that there exists a positive constant
$c(g) >0$, uniquely determined by the metric $g$ via a
Dirac-harmonic spinor, such that $m(g) \geq c(g)$ and the equality
$m(g) = c(g)$ occurs only if, outside a compact
set, $(N, g)$
is conformally flat  and the scalar curvature $S_g \equiv 0$ vanishes.
It will also be shown that the equality
$m(g) = c(g)$ is indeed attained if $(N, g)$ is conformally flat, the conformal factor being constant on the inner boundary $\partial N$, and
the scalar curvature is everywhere zero.
The idea to prove the rigidity statement is
that, near infinity,  one can conformally deform the considerd
metric as well as the connection, using the length of a
harmonic spinor without zeros as the conformal factor.

\bigskip  \noindent
\section{The Witten-Herzlich method}

In this section we recall some basic facts concerning the Witten-type spinor method used by Herzlich to prove a positive mass theorem for manifolds with inner boundary [2, 5, 6, 9, 11].
Let $( \partial_{\theta},  \partial_{\phi},  \partial_r)$ be a frame field on $(N, g)$ determined by spherical coordinates $( \theta,  \phi , r)$. Applying the Gram-Schmidt orthogonalization process to $( \partial_{\theta},  \partial_{\phi},   \partial_r)$, we obtain a $g$-orthonormal frame
$(E_1, E_2, - E_3)$, defined on an open dense subset of $N$, such that $V := - E_3$ is the outward unit normal to hypersurfaces $(S(r), g), \, r \geq r_o$,
and each $E_j, \, j=1,2$, is tangent to $S(r)$, where $( S(r), g)$ denotes hypersurface $S(r)$ equipped with the metric induced by $g$.  Let $\nabla$ and $\nabla^{\partial}$  be the Levi-Civita connection of $(N, g)$ and $(\partial N, g)$, respectively.  Let $D$ be the Dirac operator of $(N, g)$ and $D^{\partial}$ the induced Dirac operator of $( \partial N, g )$, respectively.     Let  $\Theta :=  \nabla V$ be the second fundamental form of $(\partial N, g)$. Then we have
\[   \nabla_X \psi =  \nabla^{\partial}_X  \psi + \frac{1}{2} \Theta(X) \cdot E_3 \cdot \psi   \]
for all vectors $X$ on $\partial N$ and so
\begin{equation}
 D \psi   -  E_3 \cdot
\nabla_{E_3}  \psi                     =
 \sum_{i=1}^{2}  E_i   \cdot  \nabla^{\partial}_{E_i}  \psi
- \frac{1}{2}  ( {\rm Tr}_g \Theta )  E_3 \cdot \psi  .
\end{equation}

\bigskip   \noindent
Let $\Sigma(N)$ and $\Sigma( \partial N )$ be the spinor bundle of $(N, g)$ and $( \partial N, g )$, respectively.
Recall that the Clifford bundle ${\rm Cl} (\partial N) $ may be thought of as a subbundle of  ${\rm Cl} (N) $,
the Clifford multiplication
$ {\rm
Cl}(\partial N)  \times  \Sigma(\partial N)   \longrightarrow \Sigma(\partial N) $ being
naturally related to the one  $ {\rm Cl}(N)  \times  \Sigma(N)
\longrightarrow \Sigma(N) $  via either
\begin{equation}
\pi_{\ast} (E_i  \cdot   E_3  \cdot  \psi)   =   E_i   \cdot  ( \pi_{\ast}
\psi ) ,   \quad  i=1,2 ,
\end{equation}      or
\begin{equation}
-  \pi_{\ast} (E_i  \cdot   E_3  \cdot  \psi)   =   E_i   \cdot  ( \pi_{\ast}
\psi ) ,
\end{equation}
where $\pi_{\ast} :  \Sigma(N)   \longrightarrow   \Sigma(\partial N)$ is
the restriction map.
The equation (2.1) is then projected to $\partial N$ as
\begin{equation}
\pi_{\ast} ( E_3 \cdot D \psi + \nabla_{E_3}  \psi  )  =   \mp   \sum_{i=1}^{2}  D^{\partial} (\pi_{\ast} \psi) +
\frac{1}{2}  ( {\rm Tr}_g {\rm II} )  (\pi_{\ast} \psi)  .
\end{equation}

\noindent
Regarding $\nabla^{\partial} \psi, \,    \psi   \in  \Gamma( \Sigma (\partial N) )$, as spinor fields on $N$, not projected to the boundary  $\partial N$,   one verifies easily that the formula
\[    \nabla^{\partial}_X (  E_3 \cdot \psi )  =  E_3    \cdot  \nabla^{\partial}_X  \psi      \]
makes sense. Therefore $D^{\partial}$ anticommutes with the action of the unit normal $E_3$, and hence  the discrete eigenvalue spectrum of
$D^{\partial}$ is symmetric with respect to zero. Moreover, we note  that, since the  smallest absolute value of eigenvalues of $D^{\partial}$ must satisfy
\begin{equation}
\lambda   \   \geq     \   2 \,   \sqrt{  \frac{\pi} {{\rm Area}(\partial N, g) }   } ,
\end{equation}
there is no non-trivial solutions to the equation $D^{\partial} \varphi = 0$.

\bigskip
Let
$( \cdot  ,  \cdot  )_g = {\rm Re} \langle  \cdot  ,   \cdot  \rangle_g$  be the real part of the standard Hermitian product $\langle  \cdot ,  \cdot \rangle_g$ on the spinor bundle $\Sigma(N)$ over $(N, g)$. Then, using the scalar product $( \cdot  ,   \cdot  ) = (  \cdot ,  \cdot  )_g$, one can
describe the asymptotic behaviour of spinor fields as
\begin{equation}
    \vert \psi  \vert  = \sqrt{ ( \psi,  \psi ) }  =  O (r^{-\kappa}),  \quad    \vert  \nabla  \psi  \vert
=  O (r^{-1-\kappa}) ,   \quad    \mbox{etc.},
\quad  \kappa >0  .
\end{equation}

\bigskip  \noindent
{\bf  Remark:}   Using the formulas in Proposition 2.1 and Proposition 2.3 of the paper [8],
one verifies that (2.6) is in fact equivalent to the decay condition
\[   \vert  \psi \vert_{g_{\rm eu}} = \sqrt{ ( \psi,  \psi )_{g_{\rm eu} }  } =  O (r^{-\kappa}) ,   \quad
 \vert    \nabla^{g_{\rm eu}} \,  \psi  \vert_{g_{\rm eu}}
=  O (r^{-1-\kappa}) ,    \quad       \mbox{etc.}  ,     \]
described in terms of the flat metric $g_{\rm eu}$.

\bigskip
Let $P_{\pm}$
be the $L^2$-orthogonal projection onto the subspace of positive (resp. negative) eigenspinors of the induced Dirac operator $D^{\partial}$.
Let   $W^{1, 2}_{-\tau}$ be the weighted Sobolev space
defined in [2].
In the rest of the paper, we fix a constant spinor  $\psi_o$ with $\vert  \psi_o  \vert = 1$ (i.e., $\psi_o$ is a parallel spinor with respect
to the flat metric $g_{\rm eu}$), all the components of which are constant with respect to a spinor frame field induced by rectangular coordinates,
and we use the rule (2.2) for the Clifford multiplication.
Now we
consider the PDE system :
\begin{equation}
   D \psi = 0,    \quad   \mbox{with boundary condition}  \quad  \lim_{\vert x \vert \rightarrow \infty}  \psi (x) =  \psi_o ,     \quad  P_- \psi = 0 ,
\end{equation}
where $\psi $ is a section of  $\Sigma (N)$ with
 $\psi - \psi_o    \in   W^{1, 2}_{-\tau},   \  \tau  >  \frac{1}{2} .$
(If one uses the rule (2.3) for the Clifford multiplication, then
 the spectral boundary condition $P_- \psi = 0$ must be replaced by $P_+ \psi = 0$ to gurantee positivity of the boundary term in the equation (2.8) below for the mass).

\bigskip \noindent
\begin{pro} {(see [5])} Let $(N, g)$ be a  Riemannian asymptotically flat 3-manifold of order $\tau  >  \frac{1}{2}$.
Let the scalar curvature $S_g$ of $(N, g)$ be non-negative and the mean curvature ${\rm Tr}_g(\Theta)$ of the boundary $(\partial N, g)$ satisfy
\[    \lambda   \   \geq    \     \frac{1}{2}   \,  \sup_{\partial N}  \{ {\rm  Tr}_g (\Theta)  \}  ,    \]
where $\lambda$ is the smallest absolute value of eigenvalues of the induced Dirac operator $D^{\partial}$. Then there exists a unique solution to
the PDE system (2.7).
\end{pro}

\bigskip  \noindent
Let $\psi$ be a solution to the system (2.7).   Let  $\mu_{S(r)}(g), \, \mu_{\partial N}(g), \, \mu_N(g)$ denote the volume form of $(S(r), g), \,
(\partial N, g),  \, (N, g)$, respectively. Then, applying Stokes' theorem,
 the Schr\"{o}dinger-Lichnerowicz formula and the spectral boundary condition, we have
\begin{eqnarray}
 m(g)  & =   &    \frac{1}{8 \pi}  \, \lim_{r
\to \infty} \int _{S(r)}  g ( {\rm grad}_g (\psi, \psi),  \, V ) \mu_{S(r)} (g)      \nonumber     \\
&       &       \nonumber     \\
& = &   \frac{1}{4 \pi}  \,  \int _{\partial N}  \Big(  D^{\partial} (\pi_{\ast}  \psi) -  \frac{1}{2}  {\rm  Tr}_g (\Theta) ( \pi_{\ast} \psi ) ,
\,    \pi_{\ast} \psi     \Big) \mu_{\partial N} (g)      \nonumber       \\
&       &     \nonumber       \\
&       &      +  \frac{1}{4 \pi}  \,  \int _N    \Big\{   ( \nabla \psi,    \nabla \psi )  +  \frac{1}{4}   S_g (\psi,  \psi)   \Big\}
\mu_N(g)      \nonumber         \\
&    \geq    &   \frac{1}{4 \pi}  \,  \int _{\partial N}  \Big\{ \lambda -  \frac{1}{2}  {\rm  Tr}_g (\Theta)   \Big\}  ( \pi_{\ast} \psi ,
\,   \pi_{\ast}   \psi    ) \mu_{\partial N} (g)  ,
\end{eqnarray}
which proves the following positive mass theorem.

\bigskip  \noindent
\begin{thm}  {(see [5])}   If $(N, g)$ is asymptotically flat of order $\tau  >  \frac{1}{2}$ with
$S_g   \geq  0$ and the mean curvature ${\rm Tr}_g(\Theta)$  satisfies
\[    2  \sqrt{\frac{\pi}{{\rm Area}(\partial N, g)} }   \   \geq    \  \frac{1}{2}   \sup_{\partial N}  \{ {\rm  Tr}_g (\Theta)  \}  ,  \]
then $m(g)  \geq   0$, with equality if and only if  $(N, g)$ is flat.
\end{thm}

\bigskip \noindent
Note that, if
\begin{equation}      2  \sqrt{\frac{\pi}{{\rm Area}(\partial N, g)} }   \  \geq   \   \frac{1}{2}   \sup_{\partial N}  \{ {\rm  Tr}_g (\Theta)  \}     \quad
\mbox{and}    \quad     2  \sqrt{\frac{\pi}{{\rm Area}(\partial N,
g)} }   \  \not\equiv   \    \frac{1}{2}   {\rm  Tr}_g (\Theta)
\end{equation}
on the boundary $\partial N$, then the equality
$m(g) = 0$ of Theorem 2.1 can not be attained, and hence one may find
a reasonable positive constant $c(g) > 0$ depending on the metric $g$
with $m(g)  \geq c(g)$.   In the next section, we
investigate the situation (2.9)  and  improve the rigidity statement
of Theorem 2.1.

\bigskip  \noindent
\section{Conformal change of metric using length of a spinor without zeros as the conformal factor}

We consider a conformal metric $\overline{g} = e^f g$ on $N$ with $ f  \in  W^{1,  2}_{- \tau},  \  \tau > \frac{1}{2} $.
The scalar curvatures $S_{\overline{g}}$ and  $S_g$ are related by
\begin{eqnarray}
&     &  \triangle_g ( e^{k f} )  = - ( {\rm div}_g  \circ  {\rm grad}_g  ) ( e^{k f} )     \nonumber    \\
&     &     \nonumber     \\
&   =  &   \frac{k}{2}  e^{(k +1)f} S_{\overline{g}} -  \frac{k}{2}  e^{k f}  S_g   +  \frac{k (1 - 4k)}{4}   e^{k f}  \vert  df  \vert^2_g ,
\end{eqnarray}
where $k \in  {\mathbb R}$ is an arbitrary real number,
and the mean curvatures ${\rm Tr}_{\overline{g}} (\Theta_{\overline{g}})$ and  ${\rm Tr}_g (\Theta_g)$ on the boundary $\partial N$ are related by
\begin{equation}
{\rm Tr}_{\overline{g}} (\Theta_{\overline{g}}) = e^{- \frac{f}{2} } {\rm Tr}_g (\Theta_g) - e^{- \frac{f}{2} }  df(E_3),
\end{equation}
where $E_3$ is the inward unit normal to $(\partial N, g)$. Moreover, applying (3.1) to (1.3), one verifies that the masses $m(\overline{g})$ and $m(g)$ are related
as follows:
\begin{eqnarray}
&       &  m(\overline{g}) - m(g)          \nonumber    \\
&     &      \nonumber      \\
&  =   &  \frac{1}{k} \cdot  \frac{1}{8 \pi}   \,   \int_{\partial N} g(  {\rm grad}_g ( e^{kf} ),  \, E_3  ) \mu_{\partial N} (g)
+     \frac{1}{k} \cdot  \frac{1}{8 \pi}   \,   \int_N   \triangle_g ( e^{k f} )  \mu_N (g)          \nonumber       \\
&       &     \nonumber       \\
&  =  &    \frac{1}{8 \pi}   \,   \int_{\partial N}   e^{k f}  df (E_3)   \mu_{\partial N} (g)        \nonumber       \\
&       &        \nonumber       \\
&        &      +  \frac{1}{16 \pi}   \,   \int_N    e^{k f}  \Big(    e^f  S_{\overline{g}} -  S_g   +   \frac{1- 4k}{2}    \vert df  \vert^2_g
\Big)    \mu_N (g) .
\end{eqnarray}

\bigskip
Now  let $\Sigma(N)_g$ and $\Sigma(N)_{\overline{g}}$ denote the spinor bundle of $(N, g)$ and $(N,  \overline{g})$, respectively.  Then
there are natural isomorphisms $j : T(N)  \longrightarrow T(N)$ and $j :  \Sigma(N)_g   \longrightarrow  \Sigma(N)_{\overline{g}}$
preserving the inner products of vectors and spinors as well as the Clifford multiplication
\begin{eqnarray*}
&       &   \overline{g} (jX, jY) = g (X, Y),   \quad   \langle  j \psi_1,  j \psi_2  \rangle_{\overline{g}} = \langle  \psi_1,  \psi_2  \rangle_g ,     \\
&       &     \\
&       &  (jX) \cdot (j \psi) = j( X \cdot \psi ),   \quad  X, Y  \in  \Gamma(T(N)),   \quad    \psi, \psi_1,  \psi_2   \in  \Gamma(  \Sigma(N)_g ).
\end{eqnarray*}
We fix the notation $\overline{X} := j(X)$ and $\overline{\psi} := j (\psi)$ to denote the
corresponding vector fields and spinor fields on $(N, \overline{g})$, respectively.    For
shortness we also introduce the notation $\psi_p := e^{pf} \psi, \ p \in {\mathbb R} .$  Then, one verifies that the connections
$\overline{\nabla},  \, \nabla$ and the Dirac operators
$\overline{D},  \, D$ are related as follows.

\bigskip  \noindent
\begin{pro}       \
(i)  \ $\displaystyle  \overline{\rm grad} (e^f)  = e^{- \frac{f}{2}}  \,  \overline{ {\rm grad}(e^f) } ,$

\bigskip  \noindent
(ii)  \ $\displaystyle   \overline{\nabla}_X  \overline{\psi_p}
  =     e^{pf}  \Big\{  \overline{\nabla_X  \psi} +  \frac{4p-1}{4}  e^{-f} \,
\overline{g} ( \overline{\rm grad} (e^f),  \, X)    \overline{\psi}
- \frac{1}{4}  e^{-f}  X  \cdot  \overline{\rm grad} (e^f) \cdot  \overline{\psi}   \Big\}    ,$

\bigskip  \noindent
(iii)  \  $\displaystyle   \overline{D} \, \overline{\psi_p}  =   e^{pf}   \Big\{   e^{- \frac{f}{2}}   \overline{D \psi} +
\frac{2p+1}{2}  e^{-f}  \,  \overline{\rm grad} (e^f ) \cdot  \overline{\psi}   \Big\}    .$
\end{pro}

\bigskip
Let $\varphi = \varphi_o  + \varphi_1$ be a spinor field on $(N, g)$ with  $\vert  \varphi_o  \vert = 1$ and $\varphi_1   \in   W^{1, 2}_{-\tau},   \,   \tau  >  \frac{1}{2} $.   Since $\vert  \varphi  \vert   \longrightarrow  1$ as $r  \longrightarrow  \infty$, there exists a positive constant $r_{\ast} \geq r_o $ such that $\varphi$ has no zeros in $N(r_{\ast}) :=  \bigcup_{r \geq r_{\ast}}  S(r)$.  Define a conformal metric $\overline{g}$ on $N(r_{\ast})$ by
 \[   \overline{g}  =  (\varphi, \varphi)^q g ,     \quad  q     \in  {\mathbb R}  \, .     \]
Then the connections
$\overline{\nabla},  \, \nabla$ and the Dirac operators
$\overline{D},  \, D$ are related by
\begin{eqnarray}
  \overline{\nabla}_X  \overline{\varphi_p}
 &  =  &    (\varphi, \varphi)^{pq}   \Big\{  \overline{\nabla_X  \varphi} +  \frac{q(4p-1)}{4} (\varphi, \varphi)^{-1}
\overline{g} ( \overline{\rm grad} (\varphi, \varphi),  \, X)    \overline{\varphi}          \nonumber        \\
&       &       \nonumber       \\
&        &     \qquad   - \frac{q}{4}  (\varphi, \varphi)^{-1} X  \cdot  \overline{\rm grad} (\varphi, \varphi ) \cdot  \overline{\varphi}   \Big\}   ,    \\
&       &   \nonumber     \\
\overline{D} \overline{\varphi_p}  & =  &  (\varphi, \varphi)^{pq}  \Big\{   (\varphi, \varphi)^{-\frac{q}{2}}  \overline{D \varphi} +
\frac{q(2p+1)}{2} (\varphi, \varphi)^{-1}  \overline{\rm grad} (\varphi, \varphi ) \cdot  \overline{\varphi}   \Big\}    ,
\end{eqnarray}
where  \, $\varphi_p = (\varphi, \varphi)^{pq} \varphi$.   On the other hand, we know (see [4]) that,
if $\varphi$ is an eigenspinor of $D$ on $(N(r_{\ast}), g)$, then
\begin{eqnarray}
\nabla_X \varphi  & =  & -  \frac{1}{2} (\varphi, \varphi)^{-1}
T_{\varphi} (X) \cdot  \varphi   + \frac{3}{4} (\varphi,
\varphi)^{-1}  g(  {\rm grad}(\varphi, \varphi),  X)  \varphi      \nonumber
\\    &      &    \nonumber      \\ &      &  +   \frac{1}{4} (\varphi, \varphi)^{-1} X \cdot {\rm
grad}(\varphi, \varphi) \cdot \varphi ,
\end{eqnarray}
where $T_{\varphi}$ is the energy-momentum tensor defined by
\[   T_{\varphi} (X, Y) = ( X \cdot  \nabla_Y  \varphi + Y  \cdot  \nabla_X  \varphi,   \,   \varphi )     .  \]

\bigskip  \noindent
Making use of the equations (3.4)-(3.6),  we obtain the following proposition immediately.

\bigskip \noindent
\begin{pro}
In the notations above, we have:

\bigskip \noindent
(i) If $p= - \frac{1}{2}$ and $D \varphi = 0$, then $\overline{D}  \overline{\varphi_p} = 0$.

\bigskip \noindent
(ii) If $\overline{\nabla}_X  \overline{\varphi_p}  = 0$ and $D \varphi = 0$, then
 $p = - \frac{1}{2} $ and   $  q =  1$.

\bigskip  \noindent
(iii) If  $\overline{\nabla}_X  \overline{\varphi_p}  = 0$ with $ p = - \frac{1}{2} $ and  $   q =  1$, then $D \varphi = 0$.
\end{pro}

\bigskip \noindent
We now find that, in order to improve the rigidity statement of
Theorem 2.1, the optimal parameters $p, q,$ are
\begin{equation}  p = - \frac{1}{2} ,   \quad   q =  1    .   \end{equation}

\bigskip  \noindent
For this choice of parameters, the equation (3.4) gives
\begin{eqnarray*}
&     &  (\varphi, \varphi)^2 (  \overline{\nabla}  \,  \overline{\varphi_p},  \, \overline{\nabla}  \,  \overline{\varphi_p} )     \\
&     &       \\
&  =  &  ( \nabla \varphi,  \nabla  \varphi)  +  \frac{1}{2} (\varphi, \varphi)^{-1} ( D \varphi,  {\rm grad} (\varphi, \varphi)  \cdot \varphi )
 - \frac{3}{8}
(\varphi, \varphi)^{-1}  \vert   {\rm grad} (\varphi, \varphi)    \vert^2  .
\end{eqnarray*}

\bigskip  \noindent
Applying the Schr\"{o}dinger-Lichnerowicz formula
\[   \triangle(\varphi, \varphi) = -2 ( \nabla \varphi,  \nabla \varphi ) +2 ( D^2 \varphi,  \varphi ) -  \frac{1}{2} S_g  (\varphi,  \varphi) , \]
 where $\triangle = - {\rm div} \circ {\rm grad}$,  one proves the following lemma.

\bigskip  \noindent
\begin{lem}   For the choice (3.7) of parameters, we have
\begin{eqnarray*}
&      &  \frac{1}{2} {\rm div} \{ (\varphi, \varphi)^r {\rm grad} (\varphi, \varphi) \}      \\
&      &      \\
& =  &  (\varphi, \varphi)^r    \Big\{   (\varphi, \varphi)^2 (  \overline{\nabla}  \,  \overline{\varphi_p},  \, \overline{\nabla}  \,  \overline{\varphi_p} )    +  \frac{1}{4} S_g  (\varphi,  \varphi)  -   ( D^2 \varphi,  \varphi )     -   \frac{1}{2} (\varphi, \varphi)^{-1} ( D \varphi,  {\rm grad} (\varphi, \varphi)  \cdot \varphi )     \\
&      &      \\
&      &   \qquad   +   \frac{3}{8}
(\varphi, \varphi)^{-1}  \vert   {\rm grad} (\varphi, \varphi)    \vert^2   \Big\}    +  \frac{r}{2}  (\varphi, \varphi)^{r-1}
\vert   {\rm grad} (\varphi, \varphi)    \vert^2 ,
\end{eqnarray*}
where $r \in {\mathbb R}$ is an arbitrary real number.
\end{lem}

\bigskip
Now we can prove the main result of the paper.

\bigskip  \noindent
\begin{thm}
Let $(N, g)$ be a Riemannian asymptotically flat 3-manifold of order $\tau  >  \frac{1}{2}$.
If the scalar curvature $S_g$ of $(N, g)$ is non-negative and the mean curvature ${\rm Tr}_g(\Theta)$ of $(\partial N, g)$ satisfies
\begin{equation}    2   \sqrt{\frac{\pi}{{\rm Area}(\partial N, g)} }   \  \geq   \    \frac{1}{2}  \sup_{\partial N}  \{ {\rm  Tr}_g (\Theta)  \}  ,   \qquad
     2  \sqrt{\frac{\pi}{{\rm Area}(\partial N, g)} }   \  \not\equiv   \   \frac{1}{2}   {\rm  Tr}_g (\Theta)     ,    \end{equation}
then there exists a positive constant $c(g) >0 $ uniquely
determined by the metric $g$ (as well as  a beforehand fixed constant spinor $\psi_o$)  such that    \\ (i)  $m(g)  \geq
c(g)   $ and
\\ (ii) the equality $m(g) = c(g)$  occurs only if, outside a compact set,   $g$ is conformally
flat and  the
scalar curvature $S_g \equiv 0$ vanishes.         \par
In case that $(N, g= e^{-f} g_{\rm eu})$ is conformally flat, $ f  \in  W^{1,  2}_{- \tau},  \  \tau > \frac{1}{2}$, and the conformal factor $e^{-f}$ is constant
on the boundary $\partial N$, then the equality $m(g) = c(g)$ holds.
\end{thm}

\noindent
{\bf Proof.}
Let $\psi $ be a unique solution to the PDE system (2.7).
We choose the parameter $r = - \frac{3}{4}$ in the formula of Lemma 3.1 so as to remove the terms involving $\vert   {\rm grad} (\psi, \psi)    \vert^2 $.
Then we have
\begin{eqnarray*}
 m(g)  & =   &    \frac{1}{8 \pi}  \, \lim_{r
\to \infty} \int _{S(r)} (\psi, \psi)^{ - \frac{3}{4}}  g ( {\rm
grad} (\psi, \psi),  \, V ) \mu_{S(r)} (g)     \\ &       &
\\ & = &   \frac{1}{4 \pi}  \,  \int _{S(r_{\ast})}    (\pi_{\ast} \psi,
\pi_{\ast}  \psi)^{ - \frac{3}{4}}  \Big(  D^{\partial}  (\pi_{\ast} \psi) -  \frac{1}{2}
{\rm  Tr}_g (\Theta)  (\pi_{\ast} \psi) , \,   \pi_{\ast} \psi     \Big) \mu_{S(r_{\ast})}
(g)
\\ &       &      \\ &       &      +  \frac{1}{4 \pi} \,  \int
_{N(r_{\ast})}   (\psi, \psi)^{ - \frac{3}{4}}   \Big\{ (\psi,
\psi)^2  (  \overline{\nabla}  \,  \overline{\psi_p},  \,
\overline{\nabla}  \,  \overline{\psi_p} )               +
\frac{1}{4}   S_g (\psi,  \psi)   \Big\} \mu_{N(r_{\ast})}(g)
\end{eqnarray*}
for all sufficiently large constants $r_{\ast} \geq r_o $.
On the other hand,  we know that
\begin{eqnarray*}
 m(g)  & =   &    \frac{1}{8 \pi}  \, \lim_{r
\to \infty} \int _{S(r)}  g ( {\rm grad}  (\psi, \psi),  \, V ) \mu_{S(r)} (g)     \\
&       &       \\
& = &   \frac{1}{4 \pi}  \,  \int _{\partial N}  \Big(  D^{\partial}  (\pi_{\ast} \psi) -  \frac{1}{2}  {\rm  Tr}_g (\Theta)  (\pi_{\ast} \psi) ,
\,  \pi_{\ast} \psi     \Big) \mu_{\partial N} (g)       \\
&       &      \\
&       &      +  \frac{1}{4 \pi}  \,  \int _N    \Big\{   ( \nabla \psi,    \nabla \psi )  +  \frac{1}{4}   S_g (\psi,  \psi)   \Big\}
\mu_N(g)        \\
&    >    &   \frac{1}{4 \pi}  \,  \int _{\partial N}  \Big\{   2  \sqrt{\frac{\pi}{{\rm Area}(\partial N, g)} }   -  \frac{1}{2}  {\rm  Tr}_g (\Theta)   \Big\}  ( \pi_{\ast} \psi ,
\,   \pi_{\ast} \psi    ) \mu_{\partial N} (g)    \ >  \  0 ,
\end{eqnarray*}
since $\int_N ( \nabla \psi, \nabla \psi )  > 0$ is strictly
positive. Therefore,
 there exists a positive  constant $r_{\infty} \geq r_o $ satisfying the following two conditions:
$\psi$ has no zeros in $N(r_{\infty}) = \bigcup_{r \geq r_{\infty}}  S(r)$ and
\begin{eqnarray*}
& &   \frac{1}{4 \pi}  \,  \int _{S(r_{\infty})}    (\pi_{\ast} \psi,  \pi_{\ast} \psi)^{
- \frac{3}{4}}  \Big(  D^{\partial} (\pi_{\ast}  \psi) -  \frac{1}{2}  {\rm
Tr}_g (\Theta)  (\pi_{\ast} \psi) , \,   \pi_{\ast} \psi     \Big) \mu_{S(r_{\infty})} (g)
\\ &      &      \\ &  >  & \frac{1}{4 \pi}  \,  \int _{\partial
N}  \Big\{   2  \sqrt{\frac{\pi}{{\rm Area}(\partial N, g)} }   -
\frac{1}{2}  {\rm  Tr}_g (\Theta)   \Big\}  ( \pi_{\ast} \psi , \,   \pi_{\ast} \psi
) \mu_{\partial N} (g)   \   >   \  0 .
\end{eqnarray*}
Let $r_{\rm glb} $ be the greatest lower bound of the set of all the constants $r_{\infty}$ satisfying these two conditions and define
\[   c(g)  = \frac{1}{4 \pi}  \,  \int _{S(r_{\rm glb})}    (\pi_{\ast} \psi, \pi_{\ast} \psi)^{ - \frac{3}{4}}  \Big(  D^{\partial}  (\pi_{\ast} \psi ) -  \frac{1}{2}  {\rm  Tr}_g (\Theta)  (\pi_{\ast} \psi) ,
\,   \pi_{\ast} \psi     \Big) \mu_{S(r_{\rm glb})} (g)   .   \] Then it is
clear that the statements (i) and (ii) of the theorem are true. Now it remains to prove the last statement of the theorem.
Let $\varphi = e^{\frac{f}{2}}   \psi_o$. Then Proposition 3.1 (iii) implies $D \varphi = 0$.
Furthermore,
\[     0 =  \overline{\nabla}_{\overline{E}_i}  \,  \overline{\psi_o}  =  \overline{\nabla}^{\, \partial}_{\overline{E}_i}   \,  \overline{\psi_o} +  \frac{1}{2} \Theta_{g_{\rm eu}} ( \overline{E_i}  ) \cdot  \overline{E_3} \cdot   \overline{\psi_o} =
 \overline{\nabla}^{\, \partial}_{\overline{E}_i}  \,  \overline{\psi_o}  +   \frac{1}{2 r_o}  \overline{E_i}  \cdot  \overline{E_3} \cdot   \overline{\psi_o}   ,   \qquad    i = 1, 2, \]
gives
\begin{eqnarray*}   \nabla^{\partial}_{E_i}  ( \pi_{\ast}  \varphi )
&  =  &  -  \frac{1}{2 r_o} \,  e^{\frac{f}{2}}  \, E_i  \cdot  (\pi_{\ast}  \varphi)   +  \frac{3}{4} df(E_i)  (\pi_{\ast}  \varphi) +  \frac{1}{4}  E_i   \cdot  ( \sum_{j=1}^2 df(E_j) E_j ) \cdot   (\pi_{\ast}  \varphi)      \\
&      &       \\
&  =  &   -  \frac{1}{2 r_o} \,  e^{\frac{f}{2}}  \, E_i  \cdot  (\pi_{\ast}  \varphi)  ,
\end{eqnarray*}
since the function $f$ is constant on $\partial N$. Consequently, $\varphi = e^{\frac{f}{2}}   \psi_o$ is the unique solution to the system (2.7) and
the equality $m(g) = c(g)$ holds indeed.    \hfill{QED.}

\bigskip  \noindent
{\bf Remark:}  Let $(N, g = e^{-f} g_{\rm eu} )$ be conformally flat, $ f  \in  W^{1,  2}_{- \tau},  \  \tau > \frac{1}{2} ,$ and let  the function $f$ be constant on the boundary $\partial N$.     Assume that  $S_g \geq 0$ and
the boundary condition (3.8) is satisfied.   Then the scalar curvature $S_g$  is given by (see (3.1))
\[    \triangle_g (e^{ \frac{f}{4} }) = -  \frac{1}{8} e^{ \frac{f}{4} }  S_g   \]
and so  the mass by (see (3.3))
\[   m(g) = - \frac{1}{8 \pi} \int_{\partial N} e^{ \frac{f}{4} } df(E_3) \mu_{\partial N}(g)  +
\frac{1}{16 \pi} \int_N e^{ \frac{f}{4} } S_g  \, \mu_N (g)  .  \]
Substituting the equation (3.2) into (3.8), one verifies easily that  $- df(E_3) \geq 0, \ df(E_3) \not\equiv 0 $,
and the constant $c(g)$ in Theorem 3.1 is in fact equal to
\begin{eqnarray*}   c(g) & = & - \frac{1}{8 \pi} \int_{\partial N} e^{ \frac{f}{4} } df(E_3) \mu_{\partial N}(g)   \\
&     &     \\
& = &  \frac{1}{4 \pi}  \,  \int _{\partial N}    (\pi_{\ast} \psi,  \pi_{\ast} \psi)^{ - \frac{3}{4}}  \Big(  D^{\partial} (\pi_{\ast} \psi) -  \frac{1}{2}  {\rm  Tr}_g (\Theta)  (\pi_{\ast} \psi) ,
\,  \pi_{\ast} \psi     \Big) \mu_{\partial N} (g) ,
\end{eqnarray*}
where $\psi = e^{\frac{f}{2}}   \psi_o$ is a unique solution to system (2.7). In particular, if $g$ is the spacelike Schwarzschild metric with
\[ e^{-f} = \Big( 1 + \frac{m}{2r}  \Big)^4 ,   \quad m > 0 , \]
then a direct computation, on the minimal boundary $\partial N = S(r = \frac{m}{2})$, shows that $c(g) = m$.

\bigskip  \noindent
{\bf Remark:}  It might be possible to compare the constant $c(g)$ in Theorem 3.1 with the lower bound
\[     4  \sqrt{\frac{{\rm Area}(\partial N, g)}{\pi} }     \]
of the Penrose inequality [3, 7],
in case that the boundary $(\partial N, g)$ is minimal.
It seems that
\[   4  \sqrt{\frac{{\rm Area}(\partial N, g)}{\pi} }   \ \geq  \ c(g) ,   \]
since the boundary condition (outermost minimal surface) for the constant  $4  \sqrt{\frac{{\rm Area}(\partial N, g)}{\pi} } $
is stronger than that (minimal surface) for $c(g)$.

\bigskip \bigskip  \noindent
{\bf Acknowledgement:}  The author thanks the referee for useful suggestions. This research was supported by the BK 21
project of Seoul National University and the BK 21 project of Inha
University.

\bigskip

\begin{Literature}{xx}
\bibitem{1} R. Arnowitt, S. Deser and C. Misner, Coordinate invariance and energy expressions in general relativity, Phys. Rev. 122 (1961) 997-1006.
\bibitem{2} R. Bartnik, The mass of an asymptotically flat
manifold, comm. pure Appl. Math. 34 (1986) 661-693.
\bibitem{3} H. L. Bray, Proof of the Riemannian Penrose conjecture using the
positive mass theorem, J. Diff. Geom. 59 (2001) 177-267.
\bibitem{4} Th. Friedrich and E.C. Kim, Some remarks on the Hijazi inequality and generalizations of the Killing equation for spinors,
J. Geom. Phys. 37 (2001) 1-14.
\bibitem{5} M. Herzlich, A Penrose-like inequality for the mass of Riemannian asymptotically flat manifolds, Comm. Math. Phys. 188 (1997) 121-133.
\bibitem{6} M. Herzlich, The positive mass theorem for black holes revisited,
J. Geom. Phys. 26 (1998) 97-111.
\bibitem{7} G. Huisken and T. Ilmanen, The inverse mean curvature flow and the
Riemannian Penrose Inequality, J. Diff. Geom.  59 (2001) 353-437.
\bibitem{8} E. C. Kim, A local existence theorem for the Einstein-Dirac equation, J. Geom. Phys. 44 (2002) 376-405.
\bibitem{9} T. H. Parker and C. H. Taubes, On Witten's proof of the positive
energy theorem, Comm. Math. Phys. 84 (1982) 223-238.
\bibitem{10} R. Schoen and S.-T.Yau, Proof of the positive mass II, Comm. Math.
Phys. 79 (1981) 231-260.
\bibitem{11} E. Witten, A simple proof of the postive energy theorem, Comm.
Math. Phys. 80 (1981) 381-402.
\end{Literature}

\end{document}